\documentclass[12pt]{article}

\usepackage{amsmath}
\usepackage{amssymb}
\usepackage{exscale}
\usepackage{amsthm}
\usepackage{epsfig}
\usepackage{verbatim}
\usepackage{setspace}
\usepackage{exscale}
\usepackage{epic}
\usepackage{eepic}
\usepackage{color}
\usepackage{bbm}
\usepackage{psfrag}
\usepackage{here}

\DeclareMathOperator{\disc}{disc}

\theoremstyle{plain}
\newtheorem{theorem}{Theorem}

\newtheorem{lemma}[theorem]{Lemma}
\newtheorem{corollary}[theorem]{Corollary}
\newtheorem{remark}[theorem]{Remark}
\newtheorem{definition}[theorem]{Definition}

\newcommand{\setmid}{\,|\,}
\newcommand{\N}{{\mathbb{N}}}

\newcommand{\R}{{\mathbb{R}}}

\newcommand{\HH}{{\mathcal{H}}}
\newcommand{\KK}{{\mathcal{K}}}
\newcommand{\EE}{{\mathcal{E}}}

\newcommand{\la}{\left|}
\newcommand{\ra}{\right|}

\newcommand{\ndiv}{/ \!\!\!\!\;|\,}

\begin{document}

\title{Discrepancy of Symmetric Products of Hypergraphs}

\author{Benjamin Doerr\thanks{Max--Planck--Institut f\"ur Informatik,
  Saarbr\"ucken, Germany.}
    \and Michael Gnewuch\thanks{Institut f\"ur Informatik und
    Praktische Mathematik,
    Christian-Albrechts-Universit\"at  zu Kiel, Germany. 
    Supported by the Deutsche
    Forschungsgemeinschaft, Grant SR7/10-1.} \and Nils
    Hebbinghaus\thanks{Max--Planck--Institut f\"ur Informatik,
   Saarbr\"ucken, Germany.}}

\maketitle

\begin{abstract}
For a hypergraph $\HH = (V,\EE)$, its $d$--fold symmetric product
is $\Delta^d \HH = (V^d,\{E^d |E \in \EE\})$. We give several upper
and lower bounds for the $c$-color discrepancy of such products. In
particular, we show that the bound $\disc(\Delta^d \HH,2) \le
\disc(\HH,2)$ proven for all $d$ in [B.~Doerr, A.~Srivastav, and
P.~Wehr, Discrepancy of {C}artesian products of arithmetic
progressions, Electron. J. Combin. 11(2004), Research
Paper 5, 16 pp.] cannot be extended to more than $c =
2$ colors. In fact, for any $c$ and $d$ such that $c$ does not
divide $d!$, there are hypergraphs having arbitrary large
discrepancy and $\disc(\Delta^d \HH,c) =
\Omega_d(\disc(\HH,c)^d)$. Apart from constant factors (depending
on $c$ and $d$), in these cases the symmetric product behaves no
better than the general direct product $\HH^d$, which satisfies
$\disc(\HH^d,c) = O_{c,d}(\disc(\HH,c)^d)$.
\end{abstract}

%\maketitle

\section{Introduction}

We investigate the discrepancy of certain products of hypergraphs.
In~\cite{wiraps}, Srivastav, Wehr and the first author noted the following. For a
hypergraph $\HH = (V, \EE)$ define the $d$--fold direct product and the
$d$--fold symmetric product by
\begin{eqnarray*}
  \HH^d &:=& (V^d, \{E_1 \times \dots \times E_d \setmid E_i \in \EE\}),\\
  \Delta^d \HH &:=& (V^d, \{E^d \setmid E \in \EE\}).
\end{eqnarray*}
Then for the (two-color) discrepancy \[\disc(\HH) := \min_{\chi: V \to \{-1,1\}}
\max_{E \in \EE} \left|\sum_{v \in E} \chi(v)\right|,\] we have
\begin{eqnarray*}
  \disc(\HH^d) &\le& \disc(\HH)^d,\\
  \disc(\Delta^d\HH) &\le& \disc(\HH).
\end{eqnarray*}

In this paper, we show that the situation is more complicated for
dis\-cre\-pan\-cies in more than two colors. In particular, it depends
highly on the dimension $d$ and the number of colors, whether the
discrepancy of symmetric products is more like the discrepancy of the
original hypergraph or the $d$-th power thereof. Let us make this precise:

Let $\HH = (V, \EE)$ be a \emph{hypergraph}, that is, $V$ is some finite set and
$\EE \subseteq 2^V$. 
%M Kleine Aenderungen in den naechsten Zeilen
Without loss of generality, we will assume  that $V = [n]$ for some $n\in
\N$. 
Here and in the following we use the 
shorthand $[r] :=
\{n \in \N \setmid n \le r\}$ for any $r \in \R$.
The elements of $V$ are called {\em vertices}, those of $\EE$
{\em (hyper)edges}. For  $c \in \N_{\ge 2}$, a $c$--{\em coloring}\/ of
$\HH$ is a mapping 
$\chi : V \to [c]$. The discrepancy problem asks
for balanced colorings of hypergraphs in the sense that each hyperedge shall
contain the same number of vertices in each color. The {\em discrepancy of}
$\chi$ and the $c$--{\em color discrepancy of} $\HH$ are defined by
\begin{eqnarray*}
  \disc(\HH,\chi) &:=& \max_{E \in \EE} \max_{i \in [c]} \la |\chi^{-1}(i) \cap
  E| - \tfrac 1c |E|\ra,\\
  \disc(\HH,c) &:=& \min_{\chi : V \to [c]} \disc(\HH,\chi).
\end{eqnarray*}

These notions were introduced in~\cite{wirmcol} extending the
discrepancy problem for hypergraphs to arbitrary numbers of colors
(see, e.g., the survey of Beck and S\'{o}s~\cite{bs95}). Note that $\disc(\HH) = 2
\disc(\HH,2)$ holds for all $\HH$. In this more general setting, the product
bound proven in~\cite{wiraps} is 
\begin{equation}
  \label{eq:prod}
  \disc(\HH^d,c) \le c^{d-1} \disc(\HH,c)^d.
\end{equation}
However, as we show in this paper the relation 
$\disc(\Delta^d \HH,c) = O(\disc(\HH,c))$ does not hold in
general. In Section~\ref{seclargedisc}, we give a
characterization of those values of $c$ and $d$, for which it 
is satisfied for every hypergraph $\HH$. In particular, we
present for all $c, d, k$ such that $c$ does not divide
$d!$ a hypergraph $\HH$ having $\disc(\HH,c) \ge k$ and
$\disc(\Delta^d \HH,c) = \Omega_{d}(k^d)$. In the light
of~(\ref{eq:prod}), this is largest possible apart from   factors
depending on $c$ and $d$ only.

On the other hand, there are further situations where this worst case does not
occur. We prove some in Section~\ref{secupper}, but the complete
picture seems to be complicated.

\section{Coloring Simplices}\label{seclargedisc}

To get some intuition of what we do in the remainder, let us regard
some small examples first. For $c=2$ colors and dimension $d=2$, 
%M a solution is easily seen. - Ist nicht klar, was gel\"ost werden
%M soll.
it is easy to see that $\disc(\Delta^d \HH,c) \leq \disc(\HH,c)$
holds for arbitrary hypergraphs $\HH = (V, \EE)$. As mentioned above,
we assume for simplicity that $V = [n]$.  Now coloring the vertices
above the diagonal in one color, the ones below in the other, and
those on the diagonal according to an optimal coloring for the
one-dimensional case does the job. More formally, let $\chi: V \to
[2]$. Let $\tilde \chi: V^2 \to [2]$ such that $\tilde \chi((x,y)) =
1$, if $x < y$, $\tilde \chi((x,y)) = 2$, if $x > y$, and $\tilde
\chi((x,y)) = \chi(x)$, if $x = y$. Then $\disc(\Delta^2 \HH, \tilde
\chi) = \disc(\HH,\chi)$.  Hence $\disc(\Delta^2 \HH,2) \le
\disc(\HH,2)$. This argument can be extended to arbitrary dimension to
show $\disc(\Delta^d \HH,2) \le \disc(\HH,2)$ for all $d \in \N$.

Things become more interesting if we do not restrict ourselves to $2$ colors.
For example, it is not clear how to extend the simple above/below
diagonal approach to $3$ colors (in two dimensions). In fact, as we
will show in the following, such bounds do not exist for many pairs
$(c,d)$, including $(3,2)$. However, in three dimensions
$\disc(\Delta^3 \HH,3) \le \disc(\HH,3)$ follows similarly to the
$(2,2)$ proof above. 
Indeed, for $c=2$ and $d=2$ we divided the product set $V^2$ into the sets above 
and below the diagonal, which we want to call two-dimensional
simplices of $V^2$, and the diagonal, a one-dimensional simplex of
$V^2$. For $c=3$ and $d=3$ we divide $V^3$ into the six three-dimensional
simplices in $V^3$ that we obtain from the set 
$\{x\in V^3 \setmid x_1 <x_2 <x_3\}$ by permuting coordinates, the
six two-dimensional simplices in $V^3$ that we obtain from
$\{x\in V^3 \setmid x_1 = x_2 <x_3\}$ by permuting coordinates and 
possibly changing $<$ to $>$, and finally the one-dimensional simplex 
$\{x\in V^3 \setmid x_1 = x_2 =x_3\}$.
Now with each color we color exactly two three-dimensional and two
two-dimensional simplices of $V^3$. The vertices of the diagonal
will again be colored according to an optimal coloring for the 
one-dimensional case.

%M Der obige Absatz besitzt meiner Meinung nach mehr
%M Kausalzusammenhang mit dem voranstehenden Absatz, als 
%M der nachfolgende Absatz, den ich wegkommentiert habe. 
%To analyse the problem in arbitrary dimension, we need the notation of
%simplices in higher-dimensional grids. Let $T \subseteq \N$. In $T^2$,
%we have three simplices: The elements above resp.\ below the diagonal
%and the diagonal itself. In $T^3$, we have the simplices
%$\{(x_1,x_2,x_3) \in T^3 \setmid x_{\sigma(1)} < x_{\sigma(2)} <
%x_{\sigma(3)}\}$, where $\sigma$ is some permutation of $[3]$, and the
%(six) simplices $\{(x_1,x_2,x_3) \in T^3 \setmid x_{\sigma(1)} = x_{\sigma(2)}
%\diamond x_{\sigma(3)}\}$, where $\sigma$ is some permutation of $[3]$
%and $\diamond$ is either $<$ or $>$, and finally the diagonal
%$\{(x_1,x_2,x_3) \in T^3 \setmid x_1 = x_2 = x_3\}$. For arbitrary dimension,
%the definition is as follows.

We shall now give a formal definition of $l$-dimensional simplices in
arbitrary dimensions. A set $\{x_1,\hdots,x_k\}$ of integers with $x_1
<\hdots< x_k$ is denoted by $\{x_1,\hdots,x_k\}_<$. For a set $S$ we
put
\begin{equation*}
{S \choose k} = \{T\subseteq S\, | \,|T|=k\}\,.
\end{equation*}
Furthermore, let $S_k$ be the symmetric group on $[k]$. For $l, d\in\N$ with 
$l\leq d$ let $P_l(d)$ be the set of all partitions of $[d]$ into $l$ 
non-empty subsets. 
%$J = \{J_1,\hdots,J_l\}_<$
%of $[d]$ with $\emptyset \neq J_1 <\hdots< J_l$, where $J_\mu < J_\nu$
%if and only if $\min J_\mu < \min J_\nu$.  
Let $e_1 = (1,0,\hdots,0)$, $\hdots$,
$e_d = (0,\hdots,0,1)$ be the standard basis of $\R^d$.
For  $c\in\N$ and $\lambda \in\N_0$ we write $c\,|\,\lambda$ if there exists 
an $m\in\N_0$ with $mc = \lambda$.

\begin{definition}\label{def1}
Let $d\in\N$, $l\in [d]$ and $T\subseteq\N$ finite. 
For $J = \{J_1,\hdots,J_l\} \in P_l(d)$ with $\min J_1 <\hdots< \min J_l$
put $f_i = f_i(J) = \sum_{j\in J_i} e_j$, $i=1,\hdots,l$. 
Let
$\sigma \in S_l$. We call  
\begin{equation*}
S^\sigma_J(T) := \Big\{ \sum^l_{i=1} \alpha_{\sigma(i)} f_i(J) \,|\, 
\{\alpha_1,\hdots,\alpha_l\}_< \subseteq T \Big\}
\end{equation*}
an \emph{$l$-dimensional simplex in $T^d$}. 
If $l = d$, we simply write $S^\sigma(T)$ instead of $S^\sigma_J(T)$ (as $|P_d(d)|=1$).
\end{definition}

%M Die naechsten beiden Saetze fand ich nicht praezise genug.
Clearly, the simplices in a $d$-dimensional grid $T^d$ form a
partition of $T^d$. 
The next remark shows that the numbers of $l$-dimensional
simplices are well-understood.

\begin{remark}
If $S(d,l)$, $d,l\in\N$, denote the {\rm Stirling numbers of the
  second kind}, then
$|P_l(d)| = S(d,l)$ (see, e.g. \cite{Rio}). We have
\begin{equation}
S(d,l) = \sum^l_{j=0} \frac{(-1)^j(l-j)^d}{j!\,(l-j)!}\,.
\end{equation}
Let $T\subseteq\N$ finite. Furthermore, let $I, J\in P_l(d)$ and 
$\sigma, \tau\in S_l$. If $|T|\geq l$, 
we have $S^\sigma_I(T) \neq 
S^\tau_J(T)$ as long as $I\neq J$ or $\sigma\neq\tau$. Thus the number of $l$-dimensional simplices
in $T^d$ is $l!\,S(d,l)$.
If $|T| < l$, then there exists obviously no non-empty $l$-dimensional simplex
in $T^d$.
\end{remark}

We are now able to prove the main result of this paper.

\begin{theorem}
  \label{Prop3}
  Let $c, d\in\N$.
  \begin{enumerate}
  \item If $c\,|\,k!\, S(d,k)$ for all $k\in \{2,\hdots,d\}$, then
    every hypergraph $\HH$ satisfies
    \begin{equation}
      \label{upperbound}
      \disc(\Delta^d\HH, c) \leq \disc(\HH, c)\,.
    \end{equation}
  \item If $c\ndiv k!\, S(d,k)$ for some $k\in\{2,\hdots,d\}$, then
    there exists a hypergraph $\KK$ such that
    \begin{equation}
      \label{lowerbound}
      \disc(\Delta^d\KK, c) \geq \frac{1}{3\,k!} 
      \disc(\KK,c)^k\,,
    \end{equation}
    and $\KK$ can be chosen to have arbitrary large discrepancy
    $\disc(\KK,c)$.
  \end{enumerate}
\end{theorem}

Before proving the theorem, we state some consequences. In
particular,~(\ref{upperbound}) holds never for $c=4$. For $c=3$, it
holds exactly if $d$ is odd.

\begin{corollary}
(a) Let $d\geq 3$ be an odd number. Then $\disc(\Delta^d\HH, 3) \leq
\disc(\HH, 3)$ holds for any hypergraph $\HH$.

(b) Let $d\geq 2$ be an even number and $c= 3l$, $l\in\N$. There
exists a hypergraph $\HH$ with arbitrary large discrepancy that 
satisfies $\disc(\Delta^d\HH,c) \geq \frac{1}{6} \disc(\HH,c)^2$.
\end{corollary}

\begin{proof}
Obviously $3 |\,k!$ for all $k\geq 3$. Since $S(d,2) = 2^{d-1}-1$, 
we have $3 |\, S(d,2)$ if and only if $d$ is odd. Indeed, $2^{3-1}-1
= 3$, $2^{4-1}-1 = 7$ and if $d= k+2$, then 
$2^{d-1}-1 = 4(2^{k-1}-1) + 3$, hence $3 |\, (2^{d-1}-1)$ if and only if
$3 |\, (2^{k-1}-1)$. Hence Theorem~\ref{Prop3} proves both claims.
\end{proof}

\begin{corollary}
Let $l\in\N$ and $c=4l$. For all $d\geq 2$ there exists a hypergraph
$\HH$ with arbitrary large discrepancy such that 
$\disc(\Delta^d\HH,c) \geq \frac{1}{6}\disc(\HH, c)^2$. 
\end{corollary}

\begin{proof}
As $S(d,2) = 2^{d-1}-1$ is an odd number, we have $4\ndiv 2!\,
S(d,2)$. Applying Theorem~\ref{Prop3} concludes the proof.
\end{proof}

\begin{corollary}
Let $c\geq 3$ be an odd number and $d\geq 2$. We have 
\begin{equation}
\label{odddisc}
\disc(\Delta^d\HH,c) \leq \disc(\HH,c)\hspace{2ex}\text{for all 
hypergraphs $\HH$}
\end{equation}
if and only if we have 
\begin{equation}
\label{2odddisc}
\disc(\Delta^d\HH,2c) \leq \disc(\HH,2c)\hspace{2ex}\text{for all 
hypergraphs $\HH$}\,.
\end{equation}
\end{corollary}

\begin{proof}
  According to Theorem~\ref{Prop3}, (\ref{odddisc}) is equivalent to
  the statement that $c |\,k!\,S(d,k)$ for all $k\in \{2,\hdots,d\}$.
  But, since $2 |\,k!$ for all $k\geq 2$ and $c$ is odd, this is
  equivalent to $2c |\,k!\,S(d,k)$ for all $k\in \{2,\hdots,d\}$,
  which is equivalent to~(\ref{2odddisc}).
\end{proof}

We now prove the upper bound Theorem~\ref{Prop3}(i). The main idea is
that each hyperedge of the symmetric product intersects all
$l$-dimensional simplices with same cardinality. Hence we may color
the simplices monochromatically if we can use each color equally often
for each $l \ge 2$.

\begin{proof}[Proof of Theorem~\ref{Prop3}(i)] 
  Let $c, d$ be such that $c\,|\,k!\,S(d,k)$ for all
  $k\in\{2,\hdots,d\}$.  Let $\HH = (V, \EE)$ be a hypergraph and let
  $\psi :V\rightarrow [c]$ such that $\disc(\HH,\psi)=\disc(\HH,c)$.
  For $X\subseteq V$, put $D(X)=\{(x,\hdots,x)\mid x\in X\}$. We
  define the following $c$-coloring $\chi: V^d \to [c]$. For
  $(v,\hdots,v) \in D(V)$, set $\chi(v,\hdots,v) = \psi(v)$. For the
  remaining vertices, let $\chi$ be such that all simplices are
  monochromatic, and for each $k$ there are exactly $\tfrac 1c k!
  S(d,k)$ monochromatic $k$-dimensional simplices in each color.
  
  Let $E\in\EE$ and put $R(E) := E^d\setminus D(E)$. For any
  $k\in\{2,\hdots,d\}$ and any two $k$-dimensional simplices $S, S'$
  we have $|S\cap R(E)|=|S'\cap R(E)|$. Therefore, our choice of
  $\chi$ implies $|\chi^{-1}(i)\cap R(E)|=\frac{1}{c}|R(E)|$ for all
  $i\in [c]$. Hence
\begin{eqnarray*}
\lefteqn{\max_{i\in [c]} \Big| |\chi^{-1}(i)\cap E^d| -
  \frac{|E^d|}{c} \Big|}\\ 
&=&\max_{i\in [c]} \Big| 
|\chi^{-1}(i)\cap R(E)| - \frac{|R(E)|}{c} 
+ |\chi^{-1}(i)\cap D(E)| - \frac{|D(E)|}{c} \Big|\\
&=&\max_{i\in [c]} \Big| |\chi^{-1}(i)\cap D(E)| - \frac{|D(E)|}{c} \Big|
= \max_{i\in [c]} \Big| |\psi^{-1}(i)\cap E| - \frac{|E|}{c} \Big|\,.
\end{eqnarray*}
This calculation establishes $\disc(\Delta^d\HH,c) \leq \disc(\HH,c)$.
\end{proof}

To prove the lower bound in Theorem~\ref{Prop3}, we use the following Ramsey theoretic approach.

\begin{lemma}
\label{Ramsey}
Let $c, d\in\N$. For all $m\in\N$ there exists an $n\in\N$
having the following property:
For each $c$-coloring $\chi:[n]^d \to [c]$ we find a subset 
$T\subseteq [n]$ with $|T| = m$ such that for all $l\in [d]$ each 
$l$-dimensional simplex in $T^d$ is monochromatic with respect to
$\chi$.
\end{lemma}

\begin{proof}[Proof of Lemma~\ref{Ramsey}]
The proof is based on an argument from Ramsey theory. First
we verify the statement of Lemma~\ref{Ramsey} for a fixed
simplex. Then, by induction over the number of all simplices, we prove
the complete assertion of Lemma~\ref{Ramsey}.

\emph{Claim}: For all $m\in\N$, all $l\in[d]$, all $\sigma\in S_{l}$,
and all $J\in P_{l}(d)$, there is an $n\in\N$ such that for all
$N\subseteq\N$ with $|N|=n$ and each $c$--coloring $\chi:N^{d}\to [c]$
there is a subset $T\subseteq N$ with $|T|=m$ and $S^{\sigma}_{J}(T)$
is monochromatic with respect to $\chi$.

\emph{Proof of the claim}: By Ramsey's theorem (see,
e.g. \cite{GRS}, Section 1.2), for every $l\in [d]$ there exists an $n$
such that for each $c$-coloring $\psi: {[n]\choose l} \to [c]$ there
is a subset $T$ of $[n]$ with $|T| = m$ and ${T\choose l}$ is
monochromatic with respect to $\psi$. Let
$N\subseteq\N$ with $|N|=n$. We can assume $N=[n]$ by renaming the
elements of $N$ and preserving their order. Let $\chi: [n]^{d}\to [c]$ be
an arbitrary $c$--coloring. We define $\chi_{l,\sigma,J}:{[n]\choose
l}\to [c]$ by $\chi_{l,\sigma,J}(\{x_1,\hdots,x_l\}_<) =
\chi(\sum\limits_{i=1}^{l}x_{\sigma(i)}f_{i})$, where the
$f_{i}=f_{i}(J)$ are the vectors corresponding to the partition $J$
introduced in Definition~\ref{def1}. By the Ramsey theory argument
there is a $T\subseteq N$ with $|T|=m$ and $\chi_{l,\sigma,J}$ is
constant on ${T\choose l}$. Hence, $S^{\sigma}_{J}(T)$ is
monochromatic with respect to $\chi$. This proves the claim.

Now we derive Lemma~\ref{Ramsey} from the claim. Each simplex is
uniquely determined by a pair 
$$(\sigma, J)\in\bigcup\limits_{l=1}^{d}\left(S_{l}\times P_{l}(d)\right).$$
Let $(\sigma_{i},J_{i})_{i\in[s]}$ be  an enumeration of all these
pairs. Put $n_{0}:=m$. 
We proceed by induction. Let $i\in [s]$ be such that $n_{i-1}$ is
already defined and has the property that for any $N\subseteq\N$,
$|N|=n_{i-1}$ and any coloring $\chi: N^{d}\to [c]$ there is a
$T\subseteq N$, $|T|=m$ such that for all $j\in [i-1]$,
$S_{J_{j}}^{\sigma_{j}}(T)$ is monochromatic. Using the claim, we
choose $n_{i}$ large enough such that for each $N\subseteq\N$
with $|N|=n_{i}$ and for each $c$--coloring $\varphi:
N^{d}\to [c]$ there exists a subset $T$ of $N$ with
$|T| = n_{i-1}$ and $S_{J_{i}}^{\sigma_{i}}(T)$ is
monochromatic with respect to $\varphi$. Note that there is a
$T'\subseteq T$, $|T|=m$ such that $S_{J_{j}}^{\sigma_{j}}(T')$ is
monochromatic for all $j\in [i]$. Choosing $n:=n_{s}$ proves the lemma.
\end{proof}

Related to Lemma~\ref{Ramsey} is a result of Gravier, Maffray, Renault and
Trotignon~\cite{gravier}. They have shown that for any $m\in\N$ there is an $n \in
\N$ such that any collection of $n$ different sets contains an induced subsystem on $m$
points such that one of the following holds: (a) each vertex forms a singleton,
(b) for each vertex there is a set containing all $m$ points except this one,
or
(c) by sufficiently ordering the points $p_1, \ldots, p_m$ we have that all
sets
$\{p_1, \ldots, p_\ell\}, \ell \in [m]$, are contained in the
system.\footnote{To
  be precise, the authors also have the empty set contained in cases (a) and
(c)
  and the whole set in case (b). It is obvious that by altering
  $m$ by one, one can transform one result into the other.}

In our language, this means that any $0,1$~matrix having $n$ distinct
rows contains a $m\times m$~submatrix that can be transformed through row and
column permutations into a matrix that is (a) a diagonal matrix, (b)
the inverse of a diagonal matrix, or (c) a triangular matrix.  

Hence this result is very close to the assertion of Lemma~\ref{Ramsey} for dimension $d=2$
and $c=2$ colors. It is stronger in the sense that not only monochromatic
simplices are guaranteed, but also a restriction to $3$ of the $8$ possible
color combinations for the $3$ simplices is given. Of course, this stems from the facts
that
(a) column and row permutations are allowed, (b) not a submatrix with index set
$T^2$ is provided but only one of type $S \times T$, and (c) the assumption of
having different sets ensures sufficiently many entries in both colors.

We are now in the position to prove the second part of Theorem~\ref{Prop3}.
\begin{proof}[Proof of Theorem~\ref{Prop3}(ii)]
Let $c$ and $d$ be such that $c\ndiv k!\, S(d,k)$ for some $k\in\{2,\hdots,d\}$.
Let $m$ be large enough to satisfy 
\begin{equation*}
\frac{1}{2}{m\choose \kappa} - \sum^{\kappa-1}_{l=0} l!\, S(d,l) {m\choose l}
\geq \frac{1}{3\,k!} m^{k}
\end{equation*}
for all $\kappa \in\{k,\hdots,d\}$.
(This can obviously be done, since the left hand side of the last 
inequality is of the form $m^\kappa/2\kappa! + O(m^{\kappa-1})$ for $m\to
\infty$.)
Using Lemma~\ref{Ramsey}, we choose $n\in\N$ such that for any 
$c$-coloring $\chi: [n]^d\to [c]$ there is an $m$-point set $T\subseteq [n]$ with all simplices in $T^d$ being monochromatic
with respect to $\chi$. 

We show that $\KK=\left([n],{[n]\choose m}\right)$ satisfies our claim. Let $\chi$
be any $c$--coloring of $\KK$, choose $T$ as in Lemma~\ref{Ramsey}. Let
$\kappa\in\{k,\hdots,d\}$ be such
that for each $l\in\{\kappa+1,\hdots,d\}$ there is the same number of
$l$-dimensional simplices in $T$ in each color but not so for the
$\kappa$-dimensional simplices. With
\begin{equation*}
\mathcal{S} := \bigcup^d_{l=\kappa}\bigcup_{J\in P_l(d)}\bigcup_{\sigma\in S_l}
S^\sigma_J(T)
\end{equation*}
we obtain
\begin{eqnarray*}
\lefteqn{\disc(\Delta^{d}\KK,\chi)}\\
&\geq&\max_{i\in [c]} \Big| |\chi^{-1}(i)\cap T^d| - \frac{|T^d|}{c} \Big|\\
&\geq &\max_{i\in [c]} \bigg\{ \Big| 
|\chi^{-1}(i) \cap \mathcal{S}| - \frac{|\mathcal{S}|}{c} \Big|
- \Big| |\chi^{-1}(i) \cap (T^d\setminus \mathcal{S})| 
- \frac{|T^d\setminus \mathcal{S} |}{c} \Big| \bigg\}\\
&\geq &\max_{i\in [c]} \Big| \sum_{J\in P_\kappa(d), \sigma\in S_\kappa} 
|\chi^{-1}(i) \cap S^\sigma_J(T) | - \frac{\kappa!\,S(d,\kappa)}{c} 
{m\choose\kappa}  \Big| \\
& &- \frac{c-1}{c} \bigg( m^d - \sum^{d}_{l=\kappa} l!\, S(d,l) 
{m\choose l} \bigg) \\
&\geq &\frac{1}{2}{m\choose \kappa} - \sum^{\kappa-1}_{l=0} l!\, S(d,l) 
{m\choose l}
\geq   \frac{1}{3\,k!} m^{k}\,.
\end{eqnarray*}
This establishes $\disc(\Delta^d\KK, c) \geq \frac{1}{3\,k!} m^{k}$.
Note that our choice of $n$ implies
$\disc(\KK,c)=\left(1-\tfrac{1}{c}\right)m$. 
\end{proof}

\section{Further Upper Bounds}\label{secupper}

Besides the first part of Theorem~\ref{Prop3}, there are more ways to obtain
upper bounds.
\begin{theorem}
\label{PropPrim}
Let $\HH = (V, \EE)$ be a hypergraph. Let $p$ be a prime number, $q\in\N$ and
$c=p^q$. Furthermore, let $d\geq c$ and $s = d - (p-1)p^{q-1}$. Then
$\disc(\Delta^d\HH,c) \leq \disc(\Delta^s\HH,c)$.
\end{theorem}

\begin{corollary}\label{cor}
Let $\HH = (V, \EE)$ be a hypergraph.

(a) If $c$ is a prime number, $q\in \N$ and $d=c^q$, then 
$\disc(\Delta^d\HH, c) \leq \disc(\HH,c)$.

(b) For arbitrary $d\in\N$ there holds 
$\disc(\Delta^d\HH, 2) \leq \disc(\HH, 2)$.
\end{corollary}

Statement (a) of the corollary follows from the identity 
$c^q = 1 + (c-1)\sum^{q-1}_{j=0}c^j$ and the (repeated) use of   
Theorem~\ref{PropPrim}.
Conclusion (b) follows also from Theorem~\ref{PropPrim}. Note that
Theorem~\ref{Prop3} implies that in both 
parts of Corollary~\ref{cor} we have $c\,|\,k!\, S(d,k)$ for all
$k\in\{2,\hdots ,d\}$. Hence Corollary~\ref{cor} could also have been
proven by analysing the Stirling numbers.

\begin{proof}[Proof of Theorem~\ref{PropPrim}]
As always, we assume without loss of generality that 
$V = [n]$.
Let us define the shift operator $S: [n]^d\to [n]^d$ by
\begin{equation*}
S(x_1,\hdots,x_c,x_{c+1},\hdots,x_d) = (x_2,\hdots,x_c,x_1,x_{c+1},\hdots,x_d)\,.
\end{equation*}
It induces an equivalence relation $\sim$ on $[n]^d$ by
$x\sim y$ if and only if there exists a $k\in [c]$ with $S^kx = y$. 
Now let $x\in [n]^d$ and denote its equivalence class by $\langle x \rangle$.
Put $k = |\langle x \rangle|$. Obviously $k$ is the minimal integer in
$[c]$ 
with $S^kx = x$. 
A standard argument from elementary group theory (``group acting on a
set'') 
shows that $k\,|\,c$.
Thus either $k=c$ or $S^{p^{q-1}}x = x$. Define 
$D = \{y\in [n]^d \,|\, |\langle y \rangle| < c\}$. Then 
\begin{equation*}
\psi: D \to [n]^{s}\,,\,y\mapsto (y_1,\hdots,y_{p^{q-1}},
y_{c+1},\hdots,y_d)
\end{equation*}
is a bijection.
For a given $c$-coloring
$\chi$ of $[n]^{s}$, we define a $c$-coloring $\tilde{\chi}$ of 
$[n]^d$ in the following way: We choose a system of representatives
$R$ for $\sim$. If $x\in R$ with $|\langle x \rangle| = c$, 
we put $\tilde{\chi}(S^i x) = i$ for all $i\in [c]$. 
If $|\langle x \rangle| < c$, then
$\tilde{\chi}(y) = (\chi\circ \psi)(y)$ for all $y\in \langle x 
\rangle$. 

Let $E\in \EE$. Notice, that $x\in E^d$ implies $\langle x \rangle \subseteq E^d$,
and $x\in D$ implies $\langle x \rangle \subseteq D$. Furthermore, the 
restriction of $\psi$ to $E^d\cap D$ is a bijection onto $E^{s}$. 
Thus 
\begin{equation*}
\begin{split}
\max_{i\in [c]} \Big| |\tilde{\chi}^{-1}(i) \cap E^d| - \frac{|E^d|}{c} \Big|
\leq &\max_{i\in [c]} \Big| |\tilde{\chi}^{-1}(i) \cap (E^d \cap D)| - 
\frac{|E^d \cap D|}{c} \Big| \\
&+ \max_{i\in [c]} \Big| |\tilde{\chi}^{-1}(i) \cap (E^d \setminus D)| 
- \frac{|E^d \setminus D|}{c} \Big|\\
\leq &\max_{i\in [c]} \Big| |\chi^{-1}(i) \cap E^{s}| - \frac{|E^{s}|}{c} \Big|
+ 0\,.
\end{split}
\end{equation*}
Hence $\disc(\Delta^d\HH,c) \leq \disc(\Delta^{s}\HH,c)$. 
\end{proof}

The following is an extension of the first statement of Theorem~\ref{Prop3}.

\begin{theorem}
\label{Neu}
Let $c$, $d\in\N$, and let $d'\in\{2.\ldots,d\}$. If $c\,|\,k!\,S(d',k)$ 
for all $k\in\{2,\ldots,d'\}$, then 
\begin{equation}
\label{neuneu}
\disc(\Delta^d\HH, c) \leq \disc(\Delta^{d-d'+1}\HH,c)
\end{equation}
holds for every hypergraph $\HH$.
\end{theorem}

\begin{proof}[Proof of Theorem~\ref{Neu}]
Let $\HH=(V,\EE)$ be a hypergraph with $V=[n]$. Let $\chi:
[n]^{d-d'+1} \to [c]$ be an arbitrary $c$--coloring. We define a
$c$--coloring $\widetilde{\chi}: [n]^{d}\to [c]$. Let $z\in[n]^{d}$, $x= (z_1,\ldots,z_{d'})$, and $y=
(z_{d'+1},\ldots,z_d)$. If $z_1 = \ldots = z_{d'} =: \zeta$, put
$\widetilde{\chi}(z) =\chi(\zeta, z_{d'+1},\ldots,z_d)$. Otherwise we find $k\in \{2,\ldots,d'\}$, $J\in P_k(d')$ and 
$\sigma\in S_k$ with $x\in S^\sigma_J([n])$. Since $c\,|\,k!\,S(d',k)$, we can color
the set $D:= \{(z_{\tau(1)},\ldots,z_{\tau(d')}, y)\,|\, \tau\in S_{d'}\}$ of cardinality
$k!\,S(d',k)$ evenly by our coloring $\widetilde{\chi}: [n]^d\to[c]$.
A similar calculation as the one at the end of the proof of Theorem~\ref{PropPrim} 
establishes $\disc(\Delta^d\HH, \widetilde{\chi}) \leq \disc(\Delta^{d-d'+1}\HH,\chi)$.
\end{proof} 

\begin{remark}
The condition in Theorem~\ref{Neu} is only sufficient but not necessary for
the validity of (\ref{neuneu}), as the following example shows:

Let $c = 4$, $d \geq c$ and $d'= 3$. According to Theorem~\ref{PropPrim},
we get for each hypergraph $\HH$ that 
$\disc(\Delta^d\HH, c) \leq \disc(\Delta^{d-2}\HH,c) = \disc(\Delta^{d-d'+1}\HH,c)$.
But we have  $2!\,S(d',2) = 6 = 3!\,S(d',3)$ and $4\ndiv 6$.

This example shows also, that the methods used in the proofs of Theorem~\ref{PropPrim}
and Theorem~\ref{Neu} are different.
\end{remark}


\begin{thebibliography}{AAA}
  

  \footnotesize  
 

\bibitem{bs95}
J.~Beck and V.~T. S\'os, Discrepancy theory,
in R.~Graham, M.~Gr\"otschel, and L.~Lov\'asz, Editors, Handbook
  of Combinatorics, Elsevier, Amsterdam, The Netherlands, 1995, 1405--1446.


\bibitem{wirmcol}
B.~Doerr and A.~Srivastav,
Multi-Color Discrepancies,
Comb. Probab. Comput. 12(2003), 365-399.


\bibitem{wiraps}
B.~Doerr, A.~Srivastav, and P.~Wehr,
Discrepancy of Cartesian products of arithmetic progressions,
Electron. J. Combin. 11 (2004), Research Paper 5, 16 pp. 

  \bibitem{GRS}
    R.~L.~Graham, B.~L.~Rothschild, and J.~H.~Spencer,
    Ramsey Theory, Second Edition, Wiley, New York, USA, 1990.

  \bibitem{gravier}
    S.~Gravier, F.~Maffray, J.~Renault, and N.~Trotignon,
    Ramsey-type results on singletons, co-singletons and monotone
     sequences in large collections of sets,
    European J. Combin. 25 (2004), 719-734.


  \bibitem{Rio}
    J.~Riordan,
    An Introduction to Combinatorial Analysis, Wiley, New York, USA, 1958.

\end{thebibliography}
\end{document}